\documentclass[11pt,reqno]{amsart}
\usepackage[top=1.25in,bottom=1.25in,left=1in,right=1in]{geometry}
\usepackage{amsfonts, amsmath, amssymb, amscd, amsthm, graphicx}
\usepackage{enumitem}
\usepackage[final, 
colorlinks=true]{hyperref}
\usepackage{verbatim}
\usepackage{url}
\usepackage{hyperref}

\usepackage{tikz,tikz-cd,amscd}
\usetikzlibrary{calc}
\usepackage{tikz,mathdots,color}
\usetikzlibrary{decorations.pathreplacing,decorations.markings}
\usetikzlibrary{shapes,shapes.arrows,positioning}

\tikzset{->-/.style={decoration={
  markings,
  mark=at position #1 with {\arrow{>}}},postaction={decorate}}}
\tikzset{->-/.default=.5}

\theoremstyle{plain}

\newtheorem*{claim*}{Claim}

\theoremstyle{definition}

\numberwithin{equation}{section}
\numberwithin{theorem}{section}

\newcommand{\param}%
	{{\mathchoice{\mkern1mu\mbox{\raise2.2pt\hbox{$\centerdot$}}\mkern1mu}%
	{\mkern1mu\mbox{\raise2.2pt\hbox{$\centerdot$}}\mkern1mu}%
	{\mkern1.5mu\centerdot\mkern1.5mu}{\mkern1.5mu\centerdot\mkern1.5mu}}}

\begin{document}

\title{The grounded spirituality of mathematics }


\author{Tarik Aougab}
\address{Department of Mathematics, Haverford College, Haverford, PA 19041, USA}
\curraddr{}
\email{taougab@haverford.edu}

\begin{abstract} We propose that mathematics is done largely as a means of both coping with and celebrating our physicality in the natural world. Inspired by David Abram's \textit{Becoming Animal}, we suggest that -- in opposition to the popular image of an otherworldy Platonic realm-- the spiritual component of doing mathematics is grounded on and in the Earth. We explore this as a means of questioning the widespread adoption of AI in the workflows of research mathematics.
\end{abstract}

\date{\today}

\maketitle 

\section{The Moment}

The cute and clean summary one often hears in regards to the invention of various numerical systems is that as human civilization increased in organizational complexity, it required increasingly complex arithmetic systems of organization. The circulation of currency created the possibility of debt\footnote{I highly recommend David Graeber's \textit{Debt} to anyone who is interested in this dubious history \cite{graeber2011debt}. It especially dovetails with what we write here in that it opens by questioning the legitimacy of a certain story that economists tell one another and their students: currency naturally arose out of barter when enough people who wanted shoes but who only had wheat got sufficiently frustrated when their sole trading partner wanted cheese but only had leather. Like all good stories, it relies on a kernel of truth in that of course currency \textit{does} solve the particular problem of allowing transactions between people who may not mutually possess what the other wants. But to say that currency sprung forth \textit{in order} to solve this very problem is another thing entirely.}, necessitating the negative integers. Concrete problems associated with land survey and agriculture gave rise to plane geometry, which in turn led to the discovery of the irrationals, and so on and so forth. 

Somewhere along the line, certain pockets of the human population amassed sufficient surplus to allow for a few individuals to devote time to mathematical inventions without any grounded application in mind, and thusly we were blessed with polynomial algebra and Cartesian approaches to geometry and imaginary numbers and.... and modular forms and Hodge cycles and Abelian categories, etc. The story, while appealing, is the kind of fiction that I think only stands up to a moment of scrutiny. And yet it persists, told over again to generation after generation of young mathematicians. In my opinion, it fails on at least two levels to capture important aspects of the development of mathematics over the centuries: 

\begin{enumerate}
\item At that critical moment in the history of humanity when mathematics allegedly untethered itself from the physical realm, the story neglects to explain \textit{why} anyone still bothered to do it. 

\item It ignores a crucial spiritual dimension to the practice of mathematics present always, even before the aforementioned moment (which I will heretofore refer to as The Moment) when mathematics ceased being an exclusively practical affair. 

\end{enumerate}

We've all heard about the cult of Pythagoras and the ways in which it -- at least so goes the legend-- violently protected the primacy of rationals\footnote{Not to mention Pythagoras' other spiritual fixations, such as his theory on the transmigration of souls and his related avoidance of beans.}. Perhaps you've also heard about the apocryphal death of Archimedes, run through by a soldier when he ignored a command to stop drawing circles in the sand. One way of reconciling the ``First came the applications'' cartoon story of mathematical history with the prevalence of such tales depicting a spiritual component to mathematical practice is to simply assert that such tales must refer to eras after The Moment, once mathematicians had already been freed from the hum drum constraints of Earthly reality. Of course no one can really pinpoint the precise location in time of The Moment but it evidently occurred before Archimedes' lifetime, and so his mind was so intensely liberated from the everyday and so exclusively dedicated to exploring the abstract beauties of the Platonic realm that no one could drag it back to the ground, not even by threat of force. 

It matters less whether such stories are actually true than that they are popularly collected and believed to contain important lessons about the right sort of life for a mathematician. In any case, the basic thesis of this piece is that The Moment doesn't exist. Not only has the practice of human mathematics always been spiritual, but this spirituality has been essential for feeling into the undeniable \textit{groundedness} of humanity on the Earth. In this sense, the dichotomy between Earthly application and the cultivation of Platonic beauty is an utterly false one. I certainly don't deny that Archimedes couldn't care less about the soldier threatening him, but maybe we should pay more attention to the fact that what captured his imagination in that exchange was quite literally the sand under his feet: the way it shifts perfectly to accommodate the strokes of the writing implement, each grain pushing off the others in tandem, behaving just fluidly enough to serve as a medium for the circles in his mind's eye and just enough as a solid such that the image once created actually persists. 

To be clear, I am not claiming that Archimedes was carrying out some sort of conscious exploration of the amazing physical properties of sand and that doing geometry was just a means to that end. Rather, I am suggesting that mathematics has always been in part a form of creative expression in which humans engage with their physical, Earthly environments. It is a meditative practice that can enhance one's awareness of their connection to the natural world, and in turn, that increased sense of mindfulness can lead the mathematician to new logical insights. This is a way to address the first objection (1) above to the ``applications first'' fairy tale: even when there is no concrete real-world problem we have in mind that we are trying to solve, we do mathematics either to tap into our belonging in and on the Earth, or as a way of grieving a sense of loss of this belonging, or perhaps as a way to manage the inevitable angst that accompanies physical existence on this planet. Mathematics-- a practice that does its best to trade in absolutes and certainties-- can serve as a salve for the anxiety accompanying the finite and unpredictable life of any creature. It is a spiritual art in which humans engage to both celebrate and cope with the core truth of humanity: we belong to, are of, and are inseparable from the physical world.

\section{AI and Becoming Animal}

Before giving some examples to demonstrate what I am talking about, let me address the title of the piece and the main sources of inspiration behind it. I believe that as a group, modern mathematicians are terrible at self-reflection, if for no other reason than having never devoted the necessary time for developing that skill at a discipline-wide scale. This weakness has allowed for extremely disconnected forms of mathematics to flourish. Instead of being in touch with our own psychologies, we chronically confuse proxies and corollaries  of mathematical practice with the root reasons for its continued existence. We are all intimately familiar with these less than satisfying explanations and they span the spectrum between ``applications first'' and ``pure first'', e.g. ``\textit{we do mathematics because it is \href{https://webhomes.maths.ed.ac.uk/~v1ranick/papers/wigner.pdf}{unreasonably effective} at solving real world problems}'', or ``\textit{we do it because it is beautiful}''. This confusion, in addition to the many ways that modern life disconnects us from ourselves and especially from our non-human counterparts, has allowed a spiritually hollow form of mathematics to dominate our intellectual landscape. 

Even our most celebrated emissaries have little practice seriously articulating to themselves why they do mathematics and as a result they tend to trade in the simplest, most boring, and imperial metaphors: \textit{it is like groping around in a dark jungle, slowly increasing one's field of view with penetrating strikes of the scythe}, or \textit{it is like being the first to ascend a peak and survey all of the land below before hoisting a flag in the crest to claim it}. Of course there are outliers\footnote{Beyond outliers, there are also of course entire mathematical traditions-- usually rooted in Indigenous practices-- that avoid these monotonous pitfalls. Check out the very interesting \href{https://indigenousmathematicians.org/}{Indigenous Mathematicians page} for more.}, for example Grothendieck's description of his creative process in terms of soaking a chestnut in water until it softens and opens on its own. But (a) Grothendieck is easily the exception to the rule when it comes to mathematicians raised in the Western world, both in terms of a commitment to anti-imperialism and in terms of cultivating a conscious awareness of the natural world, and (b) even this metaphor-- as peaceful and meditative as it surely is-- only speaks to half of the reciprocal nature between the mathematician and the mathematics. Just as we impact and are impacted by the natural world that cradles us and of which we are a piece, so too does the mathematics penetrate \textit{us}, shifts the way we think, and alters out relationship to the theory and to other mathematicians.

In any case, debates about the nature of the field in the wake of the introduction of AI is beginning to alleviate the allergy that modern mathematicians have to thinking about ourselves, but if we want to weather the turbulence that this new technology brings, there is no time to lose. We need to get much cleverer and more imaginative, quickly, about what aspects of our collective psychologies are soothed by doing mathematics, because I think this is ultimately at the core of why mathematics is done. I think of this as being complementary to David Bessis' call for mathematicians to become much more transparent (to both the public and to ourselves) about the ways that we think \cite{bessis2024mathematica}. While he rightly asks us to try to be more clear about what's going on in our minds when we're working out a new theorem, the ask here is that we try to become clear about our subconscious motivations for thinking about that theorem in the first place. And if it's the case-- as I suspect it is-- that it relates to our desire to feel at one with the rest of the natural world, then it's possible that by adopting AI into our practice we might usher in an explosion of new mathematics while simultaneously divorcing ourselves even further from our very reason for doing it. 

One easily predictable objection to this line of thought is that to consider mathematics to be a largely spiritual endeavor and not one that is meant to solve real problems and improve real people's lives is to navel-gaze. But first of all, this complaint falls into the same trap of The Moment: the spiritual approach that I am suggesting \textit{embraces} resolving real-world problems because-- unlike the popular descriptions of the mathematician with their head in the clouds, thinking up new abstractions because they are beautiful for a reason that no one can quite articulate-- this is a spirituality \textit{grounded} in and on the Earth and our connections to one another. Second, I do not honestly believe that what makes the wheels of pure mathematics turn is that its practitioners are confident they are contributing to the solutions of real world problems that will positively impact human life in the future. I am interested in exploring why mathematics \textit{is} done, not why it should be done. 

One easily predictable misinterpretation of what I hope to express is a conflation of the type of spirituality I am interested in with a sort of New Ageism, that is, an individualistic pursuit loosely constructed from a poorly misunderstood tapestry of Eastern mysticism, occultism, and various Indigenous traditions. What I want to center is a spirituality that-- in contrast to this misguided and unfortunate strain of Western consumerism (because ultimately that is what it is-- a means of selling books and courses)-- is rooted in our connections to one another and to the non-human world of which we are a small piece. I would also like to make a clear distinction between this grounded spirituality and the writings of some of the early naturalists, such as many of the White environmentalists who championed the removal of indigenous populations from (what they conceived as) the wilderness, so as to preserve and protect it. Nothing could be more opposed to a way of being that prioritizes our connections to one another. 

Another major source of inspiration for this piece is an incredible book by philosopher and ecologist David Abram, ``Becoming Animal.''\cite{abram2010becoming} It's difficult to do the book justice with any sort of description, and to be honest, perhaps one of the reasons why I found it so mind-blowing is that aside from the vaguest idea I was largely ignorant of the diversity of indigenous practices for cultivating and maintaining the relationships between humans, land, animals, and plants before reading it. But it's a deeply beautiful, poetic, and intelligent meditation on consciousness and how it is perhaps born out of the web of inter-relations between the inherent sentiences of all natural things. Because I know my people and can predict the ways in which some of those words might turn them off, let me be clear that Abram's approach is by no means anti-tech or anti-modern science. In his own words, 

\vspace{2 mm}

\textit{Nothing I write in these pages is meant to disparage the elegant disclosures revealed by our sciences, or even the detached states of mind necessary to derive those insights. I mean only to point out that the detached stance proper to science is itself dependent upon a more visceral reciprocity between the human organism and the world. It is this primordial and ongoing participation that has gradually shaped the glistening eyes with which we now view the lights and shadow of this sphere, as it has shaped our listening ears, able to detect a wide range of sounds vital to human life while oblivious to many others...}
\vspace{2 mm}

Given what faces us, I strongly recommend that mathematicians give this book a try. Even without the urgency that AI represents, I couldn't push this book any harder on my colleagues because I don't think it's a coincidence that so many mathematicians seem to be drawn to nature, to the outdoors, to walking in the woods and climbing rocks and exploring ravines. Just to whet your appetite, consider for instance this brilliant excerpt in which Abram describes the \textit{sonic} component to sunrises and sunsets, reproduced in its entirety because that is how it deserves to be read (and if you wish, note also its obvious mathematical character):

\vspace{2 mm}
\textit{A sphere, suspended in the depths, its contour vague within the ink-black darkness, a black-on-black phantom so faint it may only exist in our imagination...}

\textit{But watch: slowly, a thin margin of light rounds one edge of the sphere, outlining the orb's curve against the dark like a golden bow stretched taut. The light's edge is now gliding almost imperceptibly across the disk, leaving a radiance in its wake: the bright bow slowly thickening into a crescent of luminous blue. Soon other colors (white, brown, green) shape themselves within the blue as the shining frontier advances, until fully half the orb's visage is illuminated, then three-fourths, and then the entire face gleams perfectly round within the fathomless deep. }

\textit{Although already that face feels uncentered, unbalanced. Why? Look carefully: a new margin, it seems, is now following the first, an edge of obsidian darkness sliding ever so slowly across the face, erasing its features as it spreads. And listen: a faint sound like a whisper seems to accompany the frontier of shadow as it proceeds. Can you hear it? Move in closer, listening... Yes, the leading edge of the dark is indeed an audible as well as a visible line, a whisper rolling slowly across the sphere like a wave, drowning its colors in the anonymity of night. }

\textit{When the darkness reaches the farthest bound of the sphere, another audible wave slips in at the first side and begins its own journey, illuminating the orb as it advances, releasing the blue and the swirling white. The whispering, whistling sound that it makes as it moves seems identical to the hush of the previous, darker wave-- a hypnotic, breathy music. Let's swoop in closer still, to hear...}

\textit{Ahh, now that moving frontier sounds like a many-layered rustling, an audible curtain woven from intertwined strands of rushing wind and watery tones, some rising and then falling, others dipping and lifting, ethereal. As we drop in to the sphere-- its half-shadowed face expanding to fill our gaze-- the rushing exuberance of this strange music pours over us, its tangled strands now distinguishing themselves as whirring, fluted phrases and whistling rhythms, as the buzzy trills of towhees and the syncopated quaver of song sparrows, as winter wren's fevered twitterings and meadowlark's fluid gurglings, the loopy gossip of solitaires against a woodpecker's hammerings, loud nagging of jays and nasal warbling of house finches: all composing a dense membrane of speech sweeping slowly over the land, a sonic edge carrying warmth in its wake, an audible line with sections that fall almost silent as they pass over large bodies of water, over bays and seas, and then grow louder as they approach dry terrain. The advancing line swells with plaintive cries of gulls and nasal whistles of sandpipers, avocets, and gutteral grunts of cormorants; with the singsong of curlews, goose honks, and the thrum of hummingbirds swerving between blossoms, now the seesawing lilt of robins, hoarse phrases of tanagers, a bunting's mad rhythms, here the eerily ascending arpeggios of a Swainson's thrush, there the cascading arpeggios of a canyon wren, raucous banter of ravens (punctuated with clicks and fluid boinks), and then a loon's mournful yodeling-- each utterance with its own place, its own tonal niche in the rolling wave of discourse, within the swelling speech of the dawn as it breaks and spills the morning's light over the land. And as the wave passes and moves on, the cacophonous chorus of each place quiets down, the individual singers giving themselves over to the day and its necessary doings. }

\textit{To be sure, talk is still exchanged throughout the fields and the fragmented forests; across even the treeless shantytowns and the suburban backyards. There are the necessary check-ins, gossip, even an occasional bursting into song, yet it's far more subdued and desultory than the dawn's exaltation-- no compulsion to expend such energy now when it's needed elsewhere (finding food, raising the kids, avoiding predators). Only much later, as shadows lengthen and the light thickens to yellow and honey gold on the sides of houses, does the impulse begin to rise once again in the throat, drawn forth by the slowly sinking sun. And so from the fences, and telephone wires, and the upper branches of cotton-woods, a soaring crowd of voices lights the late-afternoon air, an interlace paean of ebullience and earnest pleasure swelling with the precise approach of evening's edge, as though it were the very voice of that edge. Which, of course, it is. }

\textit{Until, as dusk dims into night, the choral abundance faces into a quietude much deeper than the muted talk of the day, a cool hammock of silence slung between the soft scratching of crickets and the far-off hoot of a great horned owl. }

\textit{Such is the strange world we inhabit: an immense sphere around whose surface two long lines of birdsong are steadily sweeping-- always opposite one another, two breaking waves of vocal exuberance rolling ceaselessly around the planet.}

\vspace{2 mm}

\section{Examples}

It is about time that I provide some concrete examples to support my claims. I will focus on a topic near and dear to my own mathematical life, hyperbolic geometry and negatively curved spaces, although no deep knowledge or background in the area will be required to appreciate the points made. 

\begin{figure} [h]
\begin{center}
\includegraphics[height=8cm]{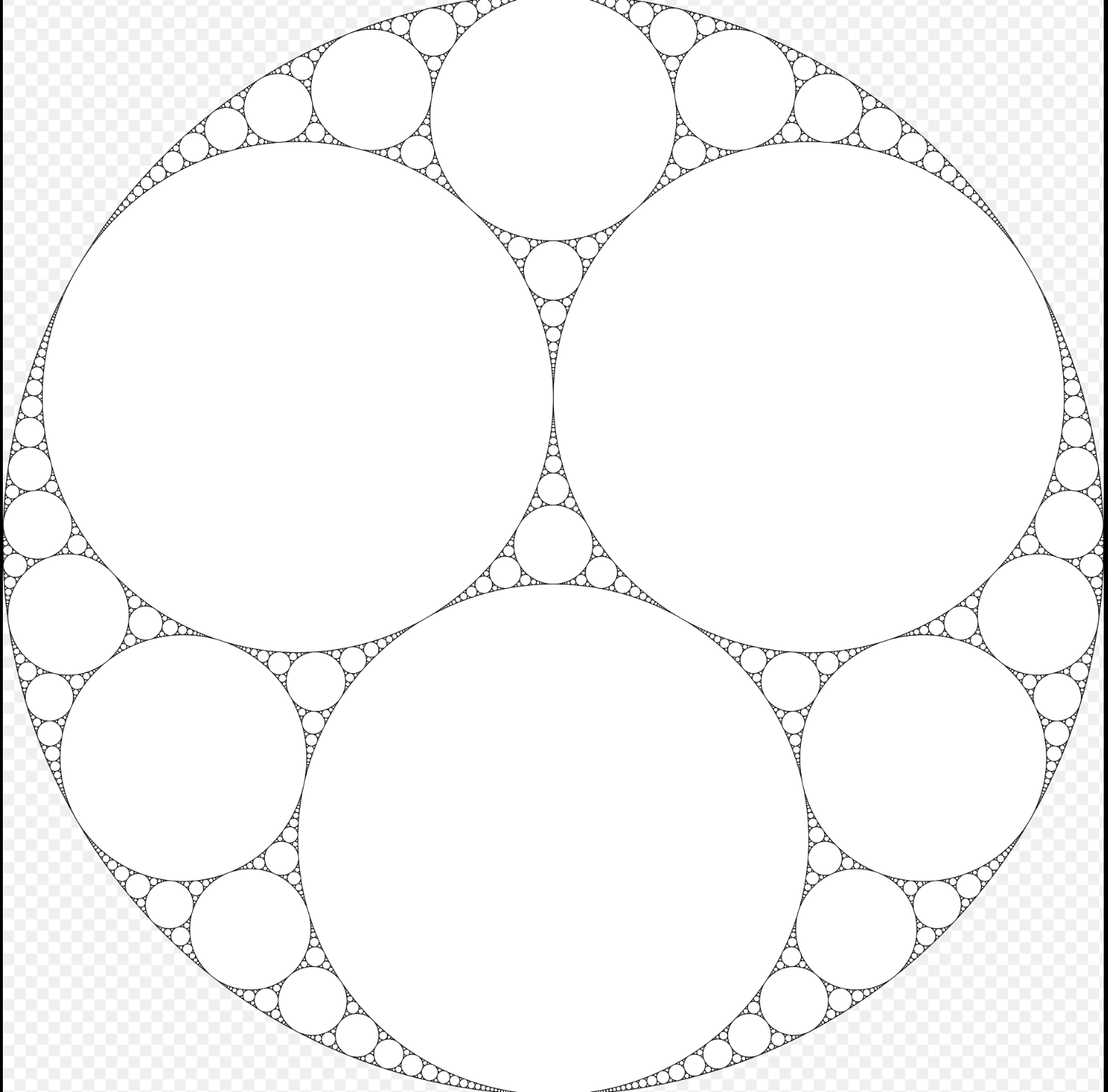}
\end{center}
\caption{The Apollonian gasket}
\label{fig:Apollonian}
\end{figure}

\begin{figure} [h]
\begin{center}
\includegraphics[height=8cm]{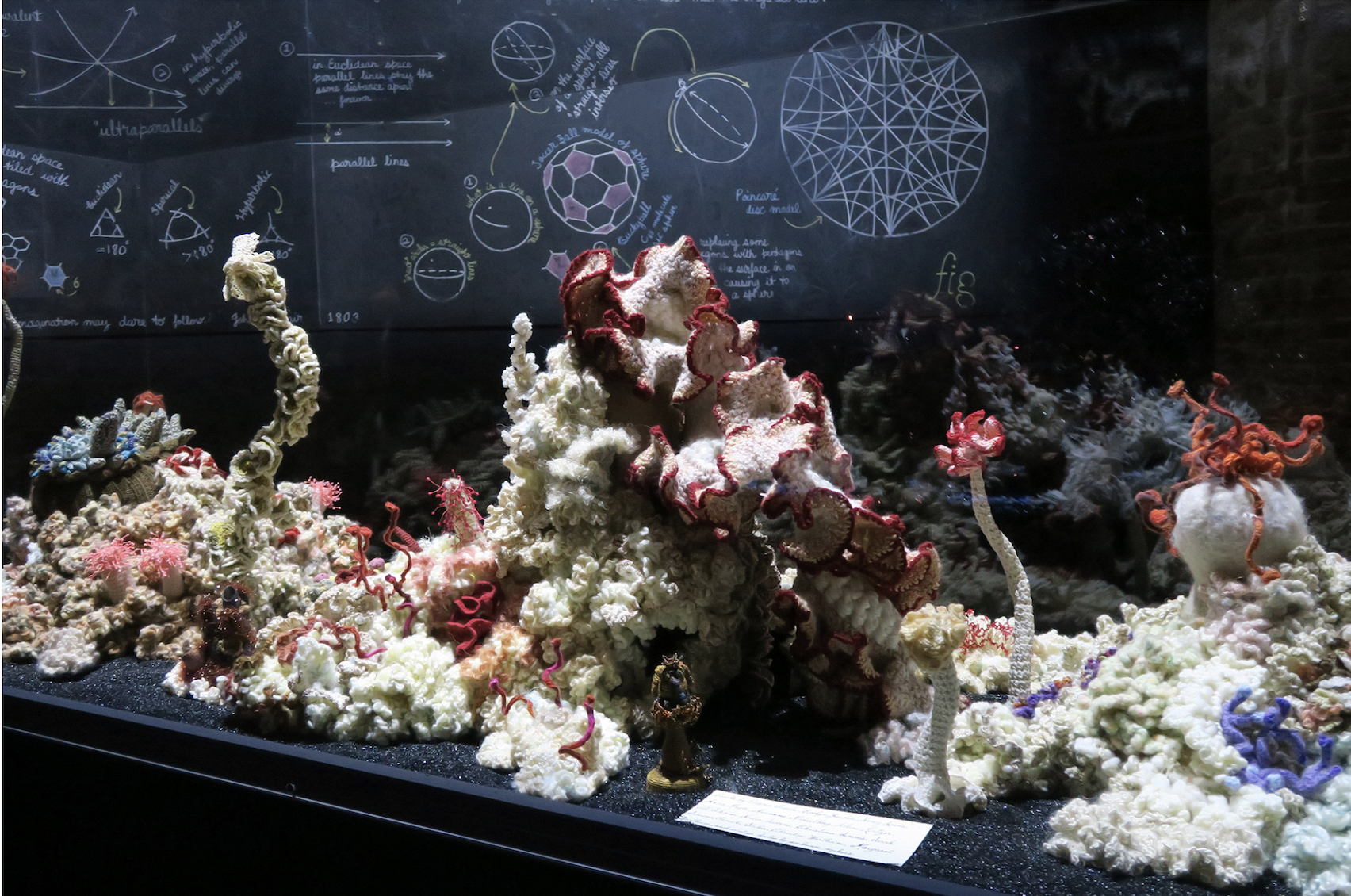}
\end{center}
\caption{A crocheted coral reef-- notice the corrugations of the coral whose boundary is marked with a deep red color.}
\label{fig:coral}
\end{figure}

There are of course many reverberations of negative curvature phenomena to be found in the natural world, for example: 

\begin{itemize}
\item The bubbles of air that form on the surface of a recently disturbed fluid, under the right conditions, may spawn to iteratively fill the space between pre-existing larger bubbles. The resulting pattern may resemble that of an Apollonian gasket, which itself relates to the action of a discrete group of isometries on hyperbolic $3$-space. If you are unfamiliar with these notions, Figure \ref{fig:Apollonian} \footnote{Picture taken from \href{https://en.wikipedia.org/wiki/Apollonian_gasket\#/media/File:Apollonian_gasket.svg} } below is all you need.

\item The larger the surface area of each leaf of a tree, the more sunlight it can absorb. But there is a natural trade-off that comes from the possibility of very large leaves casting long shadows on other leaves. One of evolution's elegant solutions was to grant certain trees with the power to develop its leafy appendages through a process of differential growth: the rate of growth is faster at the leaf's edge than at the center. This results in leaf shapes that are wavy and corrugated at the boundary because each possesses too much area to fit neatly into a single plane. Crucially, the shadow from much of this excess area is bound up within the undulating surface of the leaf itself and therefore cannot easily cast itself on to leaves located lower down on the trunk. This dynamic in which area grows too fast to be contained in a flat plane is a fundamental feature of hyperbolic geometry. This can also be observed, for some of the same biological reasons and some different ones, in coral reefs (and many mathematically minded people have been inspired enough by this to create crochet models of coral reefs, using some basic principles of hyperbolic geometry to guide the stitching-- see Figure \ref{fig:coral} \footnote{Picture taken from \url{https://crochetcoralreef.org/}}).

\end{itemize}

As incredible and beautiful as these connections are, they are not really what I want to focus on. While of course it is true that nature exhibits examples of negative curvature-- after all, we are part of nature, as are our ideas-- I don't actually think that we are driven on a psychological level to pursue hyperbolic geometry because we want to understand leaf shapes better. Instead, I think that on a deeply subconscious level, doing hyperbolic geometry speaks to a desire to connect with the infinite, inter-connected, and magical nature of our physical surroundings. I will explore this by zeroing in on three themes that are omnipresent in many (if not all!) branches of pure mathematics; so as to make the piece more concrete I will use hyperbolic geometry to anchor the investigation but these themes are hard to avoid in any mathematical discipline: 

\begin{enumerate}
\item Symmetry.
\item Engaging with the infinite. 
\item Uniqueness and rigidity. 
\end{enumerate}

In the \href{https://drive.google.com/file/d/1BUwVHXNu5-HDdOYKppYS4OR892kYb29P/view?usp=sharing}{original version of this essay}, sections that would be more easily understood by those with a strong mathematical background were color-coded in blue. Unfortunately, the ArXiv does not support highlighted text so that feature is absent in this draft. But I have done my best to make the following sections accessible to readers that are not necessarily mathematically inclined. One can also of course follow the link to the original version at the beginning of this paragraph. 

\subsection{Symmetry: the piece knows of the whole}

If you've read any popular accounts of mathematics, you've almost surely come across experts attempting to explain the beauty and importance of symmetry. And yet, there is something sort of undefinable, unreachable, or effervescent about \textit{why} symmetry is so important. On the one hand, perhaps theorems utilizing symmetry are pursued because those are the theorems we can actually prove: symmetry, when it is present, serves as a powerful tool for simplifying calculations. On the other, mathematicians tend to like symmetry because it is ``beautiful.'' 

I am not a philosopher and I know almost nothing about the theory of aesthetics. So I don't know how to articulate what it is about symmetry that many people find to be beautiful. Of course there are lots of theories that a cursory internet search pulls up:

\begin{enumerate}
\item Our brains evolved to appreciate visually symmetric objects because they require less energy to process (one only needs to comprehend half of an object for which there is a line of reflective symmetry); some studies have demonstrated a calming effect associated to that reduced neural load. 
\item Asymmetry might be evolutionarily associated with danger or ill-health. For example a fruit that appears to be ripe on one side but rotting on the other is to be completely avoided. Of course, the same is true of a fruit that is rotting all over, but the biological process of rot appears quite asymmetrical to our eyes, as opposed to the more uniform and rounded appearance characteristic of healthy fruit flesh. 
\item Many have argued over the centuries that symmetry reflects the existence of a higher power, or the presence of some sort of universal truth.  
\end{enumerate}

I am not interested in refuting any of these explanations and for all I know each is correct in its own way and sheds some light on what it is about symmetry that we like so much. I merely want to suggest an additional source for that pleasing aspect: 

\vspace{2 mm}

\textit{What if we love symmetry because it represents the capacity of a small piece to comprehend the truths of the whole? What if it speaks to our desire to bridge the seemingly unbridgeable gaps between our very different minds and bodies, and between ourselves and other dimensions of the natural world?}

\vspace{2 mm}

I feel that symmetry is beautiful to us because it serves as a visual representation of the very human longing associated to a being's ability to comprehend the whole of the world of which it is an unimaginably small piece. There is something deeply unsettling (but also deeply beautiful, as we will discuss in the next section) about the irreconcilability of the various experiences of the natural world that exist from one creature to another. To imagine structures for which the entirety is reproducible from one individual constituent soothes this unease. 

\begin{figure}
\centering
\begin{minipage}{.5\textwidth}
  \centering
  \includegraphics[width=.7\linewidth]{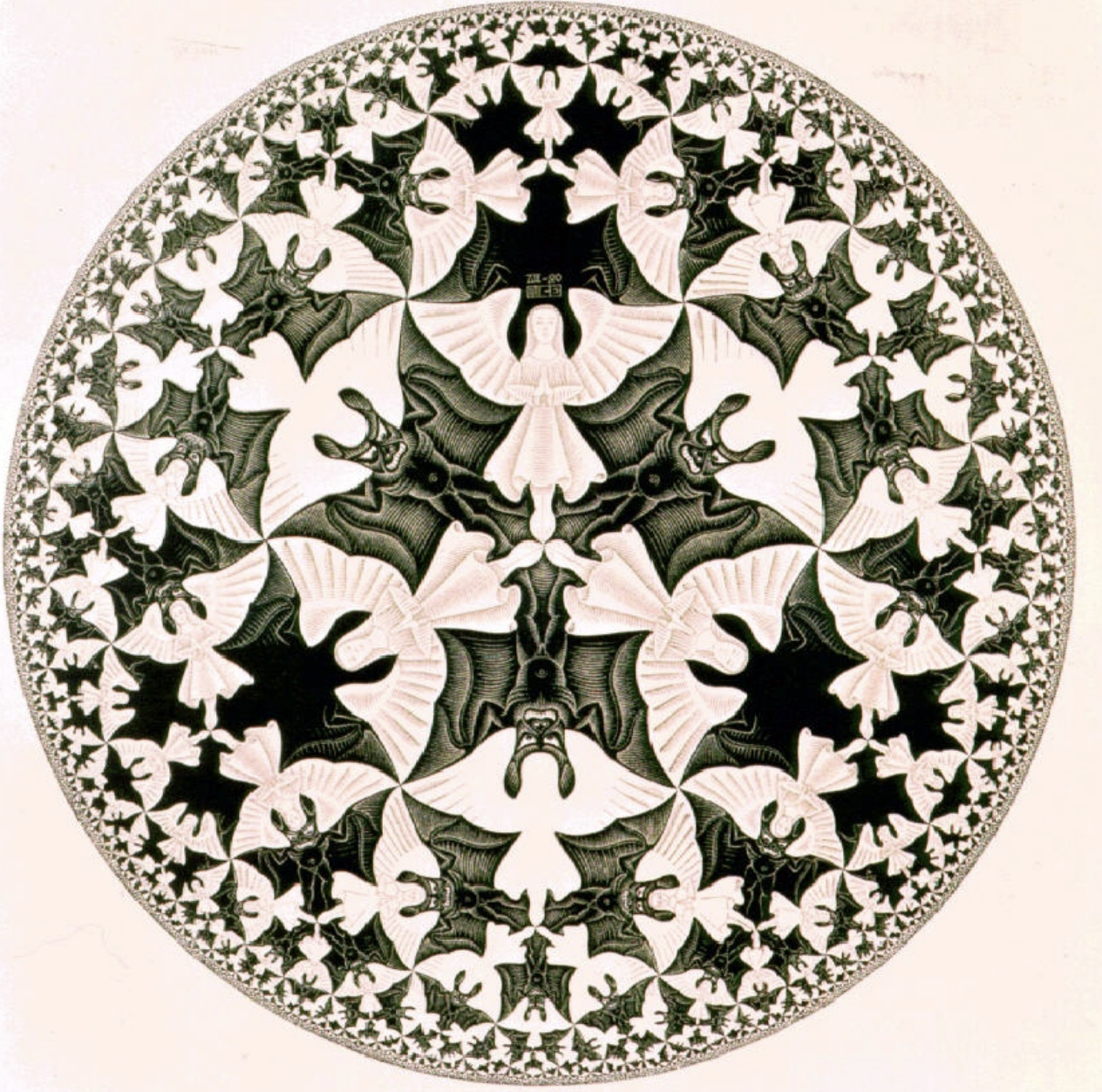}
  \caption{}
  \label{fig:Angels}
\end{minipage}%
\begin{minipage}{.5\textwidth}
  \centering
  \includegraphics[width=.7\linewidth]{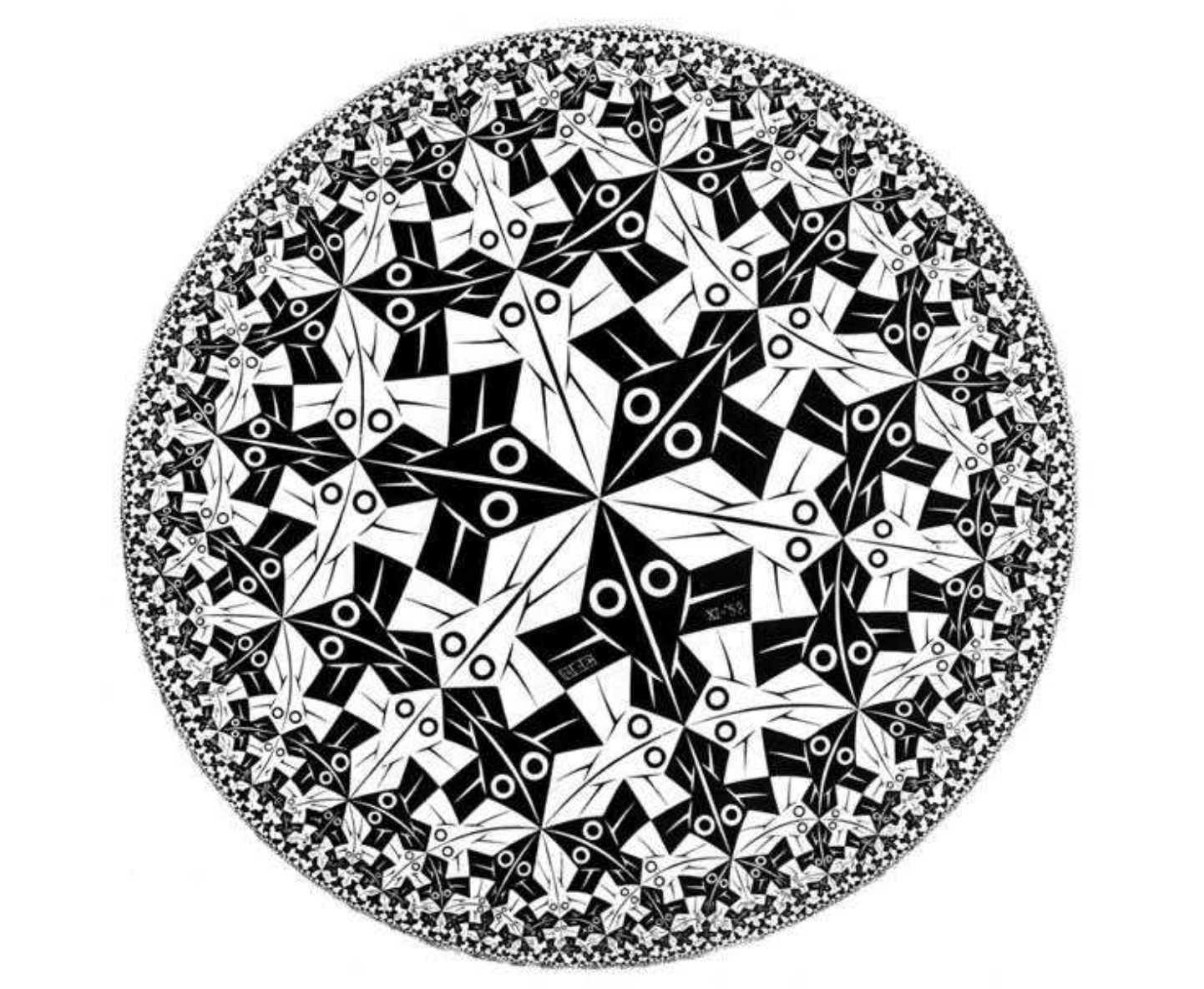}
  \caption{}
  \label{fig:fish}
\end{minipage}
\end{figure}

I can't speak for other hyperbolic geometers, but as a high school student I was inspired to learn more about the topic after coming across some of MC Escher's famed tesselations-- see Figures \ref{fig:Angels} \footnote{Figure taken from \url{https://mathoverflow.net/questions/50872/how-might-m-c-escher-have-designed-his-patterns}} and \ref{fig:fish} \footnote{Figure taken from \url{https://mathandart.com/blog/escher_and_tessellations/}} for examples. He was motivated to create them as a way to explore symmetry and order. And while he settled on the hyperbolic plane as the appropriate realm for hosting this exploration not necessarily due to considerations of symmetry\footnote{Indeed, he also created tilings of the standard flat plane, as well as works that relied on an intimate knowledge of spherical geometry.} but for reasons we will get into in the next section, hyperbolic geometry has a uniquely special relationship to symmetry.

Indeed, in two dimensions there are only three so-called \textit{Riemannian symmetric spaces}, a space that exhibits maximum isotropy in the sense that every point looks essentially identical to every other. And for the experts, all but finitely many Riemann surfaces are uniformized by the hyperbolic plane. For the non-experts, this loosely means that almost every two-dimensional expanse that is locally indistinguishable from a plane (so-called \textit{manifolds}) can be obtained by taking a chunk of the hyperbolic space and declaring a set of Pacman-like rules which dictate that a creature exiting from one edge along the chunk must re-enter seamlessly along another. In this sense, of the three symmetric spaces-- the sphere, the Euclidean plane, and the hyperbolic plane-- it is the latter that is most intimately involved in granting almost all two dimensional manifolds the opportunity to exhibit their maximal possible symmetry.

Escher's tesselations involve strikingly artistic tiles based on animals and other aesthetic objects. But when one zooms out and looks at the full tesselation, one is also naturally made aware of the various relationships \textit{between} the tiles. In truth, the description above regarding the piece knowing the truth of the whole was an overstatement, because the full tesselation can not in fact be reproduced just from one tile. This is easy to see just by thinking about the case of a single line of symmetry: if you are given half of a geometric figure and told that the full figure is obtained from it by reflecting across a line of symmetry, you can not in fact recreate the full figure unless you know where that line of symmetry is.

The identical nature of each tile naturally draws the eye to the ways in which the tiles are placed within the overall picture, relative to one another. To me, this mirrors an insight that is stamped all over the mathematics of the last 100 years: it is just as if not more important to understand the relationships between various mathematical objects as it is to understand the objects themselves. I think another aspect of this that sparks such a strong intellectual interest in us is that it dovetails with our intuitive human understanding that the spiritual truth of our world is carried not only by each individual but by the threads of kinship that sustain our collective existence.

\subsection{Infinity: the limitlessness of the physical world}

One of Abram's central invitations to the reader is to consider the possibility that every entity, no matter how seemingly small or inert, carries its own sense of experience and sentience. From this perspective, every expanse, every room in your apartment, every grove of trees is effectively boundless, for it is experienced an endless number of ways by innumerably many beings. Here is Abram describing a walk in the woods during which he feels an awareness of this truth: 

\textit{A spiderweb collapses across my face...I wipe the web from my nose and from the edge of my baseball cap; but before I can place the cap back on my head I notice the spider, so close I have to cross my eyes to see her brown body clearly; suspended by a single strand from the beak of my cap. Dangling upside down with her legs folded inward, the spider is falling away from the cap on the rapidly lengthening span. ``Hello there,'' I say hopefully, raising the baseball cap high so the spider dangles again at eye level, but she keeps riding the end of that thread down and away. I carry the cap over to a nearby scrub oak, and transfer the cap's end of the thread to an upper leaf. The veins in those leaves are clear and crisply articulated to my gaze, as are the tiny nodules toward the end of the branch where next year's buds are already preparing themselves. By holding my eyes for a few moments, the spider has drawn my senses into another scale of experience. Which opens up endlessly.}

As Abram himself notes, this frame of mind can be found represented even within the Western canon that usually so readily dismisses it. For instance, Spinoza's philosophical body of work equates God with nature itself and asserts the ensouledness of every material body\cite{Spinoza1994Ethics}. In Spinoza there is a placement of the infinite within the physical realm. This is in stark contrast with (at least the most commonly practiced and orthodox versions of) the three major religions of the Book that position God completely out of material reach; of course this had a lot to do with Spinoza's excommunication from the Portuguese synagogue. 

In any case, modern mathematics grants us a language and a set of tools for engaging logically with the infinite. We can interact and engage intellectually with it and in that sense there is something Spinoza-esque about our relationship with infinity. And to me, this has always been a deeply beautiful and alluring aspect of doing mathematics: by (literally) drawing infinity into our own world, we get to consort with boundlessness. 

As mentioned previously, characteristic to hyperbolic geometry is the exponential growth of area. One can view the ``arrival'' of infinity in finite time to the viewer's sphere of vision as a consequence of this rapid area growth. In Figures \ref{fig:Angels} and \ref{fig:fish}, one should interpret each angel/devil or fish as being of the exact same size. The ones that appear to be very small are only so because they are very far away. In that sense, the circular boundary of each figure represents a visual infinity.

This is indeed one of the key motivations behind Escher's choice of the hyperbolic plane as the stage for his tesselations; he wanted to engage directly with infinity, and doing this was made much more possible with the use of a geometry that places infinity within a bounded region to be engaged with more fully by the human mind and senses. Moreover, in hyperbolic geometry infinity itself has its \textit{own} shape and its own geometry-- in two dimensions it is the circle; in general it is a sphere of the appropriate dimension (one dimension less than the hyperbolic space it bounds). One can often glean deep truths about the dynamics of a geometric process happening within hyperbolic space by observing its imprint on the sphere at infinity. This reinforces its material physicality. 

In short, there is something irresistible about the prospect that the infinite is most naturally housed within-- as opposed to wholly outside of and orthogonal to-- our natural world. And I propose that we are drawn to mathematical traditions that, in their own ways, encapsulate this outlook.

\subsection{Rigidity: coping with impermanence}

Throughout mathematics, one finds many examples of the \textit{uniqueness} of an especially aptly defined mathematical structure. If you are careful enough about how you approach a particular mathematical property, and if the property is a sufficiently ``natural'' thing to be thinking about in the first place, you will often be able to prove that there is only one object that could possibly satisfy your requirements (of course, up to the natural notions of equivalence; we view there as being only one solution to the equation $x- 4 = 0$ even though $4$ and $3+1$ both work). 

What is it about this form of uniqueness that we find so appealing in mathematics? In hyperbolic geometry, one comes across it in many manifestations. For example, there is but one geodesic-- a path which locally minimizes the distance between points along it-- between any pair of points, or in any homotopy class of curves on a manifold modeled with hyperbolic geometry. In higher dimensions, when a closed manifold admits a hyperbolic structure, it does so in only one way. 

To me, the emotional quality of this type of uniqueness sometimes feels ethereal, frozen, or rigid. In my mind, I imagine a sort of floppy topological expanse; somewhere, in the unimaginably vast and infinitely huge universe of possible shapes that expanse can take on, there is the unique one that is of special interest because of its symmetry or its minimization of some sort of energy or length or tautness. I then ask the expanse to don this particular outfit. If it is willing, that is, if I am not too preoccupied with inclement weather or daily responsibilities or a political crisis and my mind is clear, it shimmies in and snaps into form. Once it does, it is frozen there; it can not easily deform without losing sight of the constraint that led to this crystalline avatar in the first place. 

What is pleasing about this mathematical playbook is that it flies in the face of the unpredictable and chaotic nature of everyday reality. Everything is in flux and the only constant is change itself. The moment we feel like we have a solid grasp on our surroundings is the moment it often morphs into something totally unfamiliar, new, and often scary. I think that we are naturally drawn to rigidity phenomena within mathematics for the same essential reasons that people have pondered the possible existence of permanent and unchanging spiritual realms beyond our own.

\section{Recommendations} 

If, as I propose, mathematics is done either in celebration of or in reaction to our groundedness in nature (and so in either case, mathematics revolves around this groundedness), what would that mean in terms of how we should practice it, what we should avoid, and what we should prioritize? 

The easiest and most obvious recommendation I can think of is to simply do math in nature: take a page from Archimedes and draw in the sand. While you're at it, take another page from him\footnote{Of course in reality, what we know of Archimedes was that he was engaged in extensive military scientific research, inventing a number of defensive technologies to protect Syracuse from invading Romans. There is simply no ancient equivalent to e.g. the American military which is responsible for the deaths of countless innocents all while being the single largest emitter of greenhouse gases on the planet. In any case, this sentence should be read not as a denouncement of any and all military activity, only of those militaries which act in the name of imperial domination.} and cast aspersions on the militaries of the world whose actions destroy life and the fragile equilibria upon which our Earth depends (although ideally, try not to get yourself killed while doing this). Ponder your schemes and sheaves in the woods and on the dunes and in the meadows. Hold workshops near to the wilderness and see how making this choice consciously shifts or otherwise impacts your subconscious motivations for doing mathematics. 

A stand-out example already exists: the \href{https://nesinkoyleri.org/en/main-page/nesin-maths-village/}{Nesin Math Village} in Turkey. Participants share their mathematical passions and ideas, offering courses to one another while simultaneously living on and working the land together. Since I have not gone myself I can not speak directly to the particulars of this experience, but it is something I am hoping to do in the near future. 

Wherever you are-- in your office in the middle of a big city or in the forest-- try to foster an awareness of how your mathematical practice springs forth from and depends upon your material surroundings. The air in your lungs, the fresh water you hopefully have access to, the plants and the animals and the people and the supply chains that place food in front of you and the public and private investments that fund your capacity to purchase it. Pay heed also to the important \textit{absences} that help to inspire mathematical creativity within you: if you are lucky to be doing mathematics outside of an active warzone and away from the incessant booms and clacks of bombs and guns, take that in fully. 

An essential part of doing so involves making a serious effort not to inadvertently support the very forces that steal access to those crucial absences from our mathematical siblings in Palestine, Sudan, the Congo, Cuba, Venezuela, Iran, and so many other places scarred by the insatiable death cult of imperialism. Think about adhering to the \href{https://www.justmathematicscollective.net/USACBI_statement.html}{scientific and academic boycott of Israel} if you don't already. Attend one of the monthly occurrences of the \href{https://sites.google.com/view/icm2026boycott/seminar}{Genuinely Global Mathematics Community seminar series}. Absorb the truth that mathematics is political precisely because we do it to connect ourselves to the surrounding world and to one another.

One of the reactions to the AI boom that hasn't gotten the attention it deserves is from those mathematicians who've expressed a disinterest in sitting in front of a computer all day long, no matter how much AI might bootstrap or enhance mathematical production. We should be wary of AI not merely because of its capacity to degrade our environment\footnote{AI optimists are quick to point out that the total water usage of all data centers pales in comparison to that of e.g. almond farms. But (a) almond farms are notorious water guzzlers so it's unclear that this comparison is helpful for making the point they're trying to make, and (b) this is precisely the sort of \textit{ungrounded} thinking that I am hoping to avoid, in the sense that it is not based in a particular place. One loses sight of all of the most potent political dynamics by just tallying up complete water usage for the entire world and forgetting where and from whom the majority of that water is taken. There are no almond farms in the middle of Philadelphia, but there \textit{are} talks of placing a new hyper-scale data center directly adjacent to the same predominately Black neighborhoods that have suffered from disproportionate rates of asthma and rare cancers for the last 100 years.} but also because of its potential to further detach us from that environment of which in reality we are an integral part. Even in a theoretical future that allows for the worldwide implementation of AI with zero environmental footprint, the long term impact of a technology that results in us being planted in front of screens for even longer than we already are is the severance between humans and the source of all of our creativity, passion, and joy. I don't know what the future holds but I do know that I don't want mathematics to be associated with that sort of spiritual bereftness.

As mathematicians, we are stewards of an ancient Earthly tradition in which humans meditate on the structures of pattern and shape, inspired directly by the natural world that surrounds us to filter our deepest euphorias, insights, and human anxieties through our particular medium. And for a complicated set of reasons must of us-- myself very much included-- are doing a terrible job of honoring this sacred obligation. By relegating the spiritual component of mathematics to some imagined Platonic realm, we have in essence reproduced the same banishment of spirit from its natural habitat championed by (the most prominently practiced forms of) the Abrahamic religions. And there is little I can think of that would push us faster and even further down this treacherous road than the uncritical adoption of AI. Let us reconnect to why we do what we do, and let us do so quickly.

\section{Acknowledgements}

The author sincerely thanks Jayadev Athreya, Mark Hagen, Sam Hughes, Yvonne Lai, Marissa Loving, Miriam Kuzbary, Johanna Mangahas, Maud McInerney, and Abdul Zalloum for reading a draft of this piece and offering extremely helpful feedback. 

\bibliographystyle{plainurl}
\bibliography{spiritbib}
\end{document}